\documentclass{amsart}

\usepackage{graphicx}
\usepackage[left=3cm,right=2.5cm,top=4cm,bottom=5cm]{geometry}
\newtheorem{theorem}{Theorem}[section]
\newtheorem{proposition}[theorem]{Proposition}

\theoremstyle{definition}

\theoremstyle{remark}
\newtheorem{remark}[theorem]{Remark}

\numberwithin{equation}{section}

\newcommand{\V}{\mathbb{V}}
\newcommand{\C}{\mathbb{C}}
\newcommand{\R}{\mathbb{R}}
\newcommand{\Pn}{\mathbb{P}}
\newcommand{\dS}{\mathrm{dS}}
\newcommand{\AdS}{\mathrm{AdS}}
\newcommand{\Ein}{\mathbb{E}\mathrm{in}}
\newcommand{\En}{\mathbb{E}}
\newcommand{\Hn}{\mathbb{H}}
\newcommand{\Nu}{\mathfrak{N}}

\begin{document}

\title{On the irreducible action of PSL$(2, \mathbb{R})$ on the $3$-dimensional Einstein universe. }

\author{Masoud Hassani}
\address{Universit\'e d'Avignon, Campus Jean-Henri Fabre- 301, rue Baruch de Spinoza,  BP 21239 F-84 916 AVIGNON Cedex 9.}
\email{masoud.hassani@alumni.univ-avignon.fr}
\address{University of Zanjan, Faculty of Science, Department of Mathematics, University blvd, Zanjan, Iran.}
\email{masoud.hasani@znu.ac.ir}
\thanks{This work is funded by the French Ministry of Foreign Affairs, through Campus France. 
Grant number 878286E}

\begin{abstract}
We describe the orbits of the irreducible action of PSL$(2, \mathbb{R})$ on the $3$-dimensional Einstein universe $\Ein^{1,2}$.
This work completes the study in \cite{Col}, and is one element of the classification of cohomogeneity one actions on $\Ein^{1,2}$ (\cite{Masoudthesis}).
\end{abstract}

\maketitle

\section{Introduction}

%Let $M$ be a manifold and $G$ a Lie group acting on $M$. The action of $G$ is called \textbf{cohomogeneity one} if $G$ admits a codimension one orbit in $M$.
\subsection{Einstein universe}
Let $\mathbb{R}^{2,n+1}$ denote a $(n+3)$-dimensional real vector space equipped with a non-degenerate symmetric bilinear form $\mathfrak{q}$ with signature $(2,n+1)$. The nullcone of $\mathbb{R}^{2,n+1}$ is
\begin{align*}
\Nu^{2,n+1}=\{v\in \mathbb{R}^{2,n+1}\setminus \{0\}: \mathfrak{q}(v)=0\} .
\end{align*}
The $(n+1)$-dimensional \textbf{Einstein universe} $\Ein^{1,n}$ is the image of the nullcone $\Nu^{2,n+1}$ under the projectivization:
\begin{align*}
\mathbb{P}:\mathbb{R}^{2,n+1}\setminus\{0\}\longrightarrow \mathbb{RP}^{n+2}.
\end{align*}
The degenerate metric on $\Nu^{2,n+1}$ induces a $O(2,n+1)$-invariant conformal Lorentzian structure on Einstein universe. The group of conformal transformations on $\Ein^{1,n}$ is $O(2,n+1)$ \cite{Fra}.

A lightlike geodesic in Einstein universe is a \textbf{photon}. A photon is the projectivisation of an isotropic $2$-plane in $\R^{2,n+1}$.  The set of photons through a point $p\in \Ein^{1,n}$ denoted by $L(p)$ is the \textbf{lightcone} at $p$. The complement of a lightcone $L(p)$ in Einstein universe is the \textbf{Minkowski patch} at $p$ and we denote it by $Mink(p)$. A Minkowski patch is conformally equivalent to the $(n+1)$-dimensional Minkoski space $\En^{1,n}$ \cite{Barbot}.

The complement of the Einstein universe in $\mathbb{RP}^{n+2}$ has two connected components: the $(n+2)$-dimensional Anti de-Sitter space $\AdS^{1,n+1}$ and the generalized hyperbolic space $\Hn^{2,n}$: the first (respectively the second) is the projection of the domain $\mathbb{R}^{2,n+1}$ defined by $\{ \mathfrak{q}<0 \}$ (respectively $\{ \mathfrak{q}>0 \}$).

An immersed submanifold $S$ of  $\AdS^{1,n+1}$ or $\mathbb{H}^{2,n}$ is of \textbf{signature} $(p, q, r)$ (respectively  $\Ein^{1,n}$) if the restriction of the ambient pseudo-Riemmanian metric (respectively the conformal Lorentzian metric) is of signature $(p,q,r)$, meaning that the radical has dimension $r$, and that maximal definite negative and positive subspaces have dimensions $p$ and $q$, respectively. If $S$ is nondegenerate, we forgot $r$ and simply denote its signature by $(p,q)$. 

\subsection{The irreducible representation of PSL$(2, \mathbb{R})$}
A subgroup of $O(2,n+1)$ is \textbf{irreducible} if it preserves no proper subspace of $\R^{2,n+1}$.
By \cite[Theorem 1]{Scala}, up to conjugacy, $SO_\circ(1,2)\simeq$ PSL$(2,\R)$ is the only irreducible connected Lie subgroup of $O(2,3)$.

On the other hand, for every integer $n$, it is well known that, up to isomorphism, there is only one
$n$-dimensional irreducible representation of PSL$(2,\R)$. For $n=5$, this representation is the natural action of PSL$(2, \mathbb{R})$ on the vector space $\V=\R_4[X,Y]$ of homogeneous polynomials of degree $4$ in two variables $X$ and $Y$. This action preserves the following quadratic form
$$
\mathfrak{q}(a_4X^4+a_3X^3Y+a_2X^2Y^2+a_1XY^3+a_0Y^4)=2a_4a_0-\frac{1}{2}a_1a_3+\frac{1}{6}a_2^2.
$$
The quadratic form $\mathfrak{q}$ is nondegenerate and has signature $(2,3)$. This induces an irreducible representation PSL$(2,\R)\rightarrow O(2,3)$ \cite{Col}.

\begin{theorem}\label{thm1}
The irreducible action of PSL$(2,\R)$ on the $3$-dimensional Einstein universe $\Ein^{1,2}$ admits three orbits:
\begin{itemize}
  \item An $1$-dimensional lightlike orbit, i.e. of signature $(0,0, 1)$
  \item A $2$-dimensional orbit of signature $(0, 1, 1)$,
  \item An open orbit (hence of signature $(1,2)$) on which the action is free.
\end{itemize}
\end{theorem}

\begin{figure}[h]
\centering
\hspace{-3.5mm}
\includegraphics[scale=0.18]{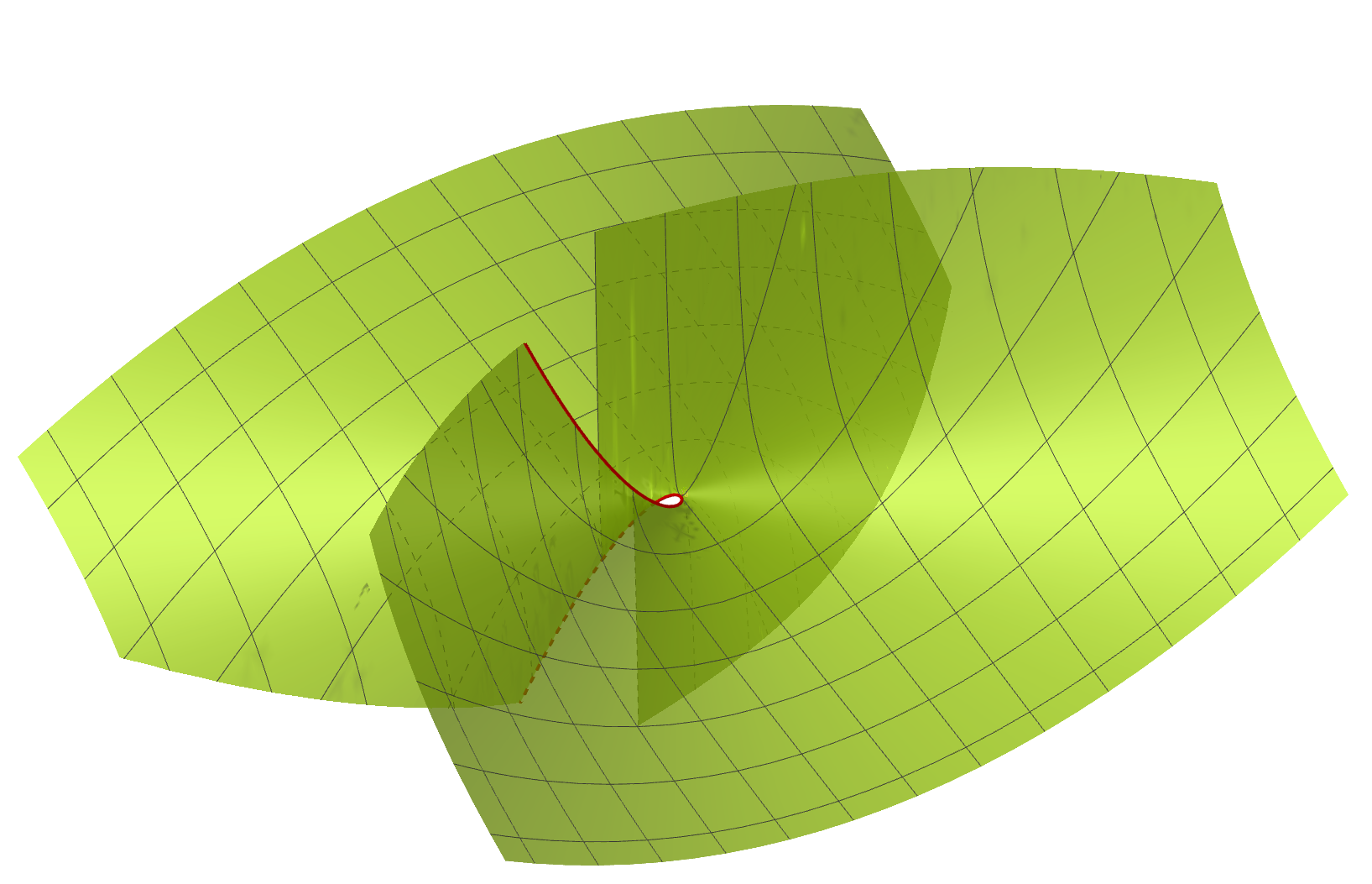}
\includegraphics[scale=0.18]{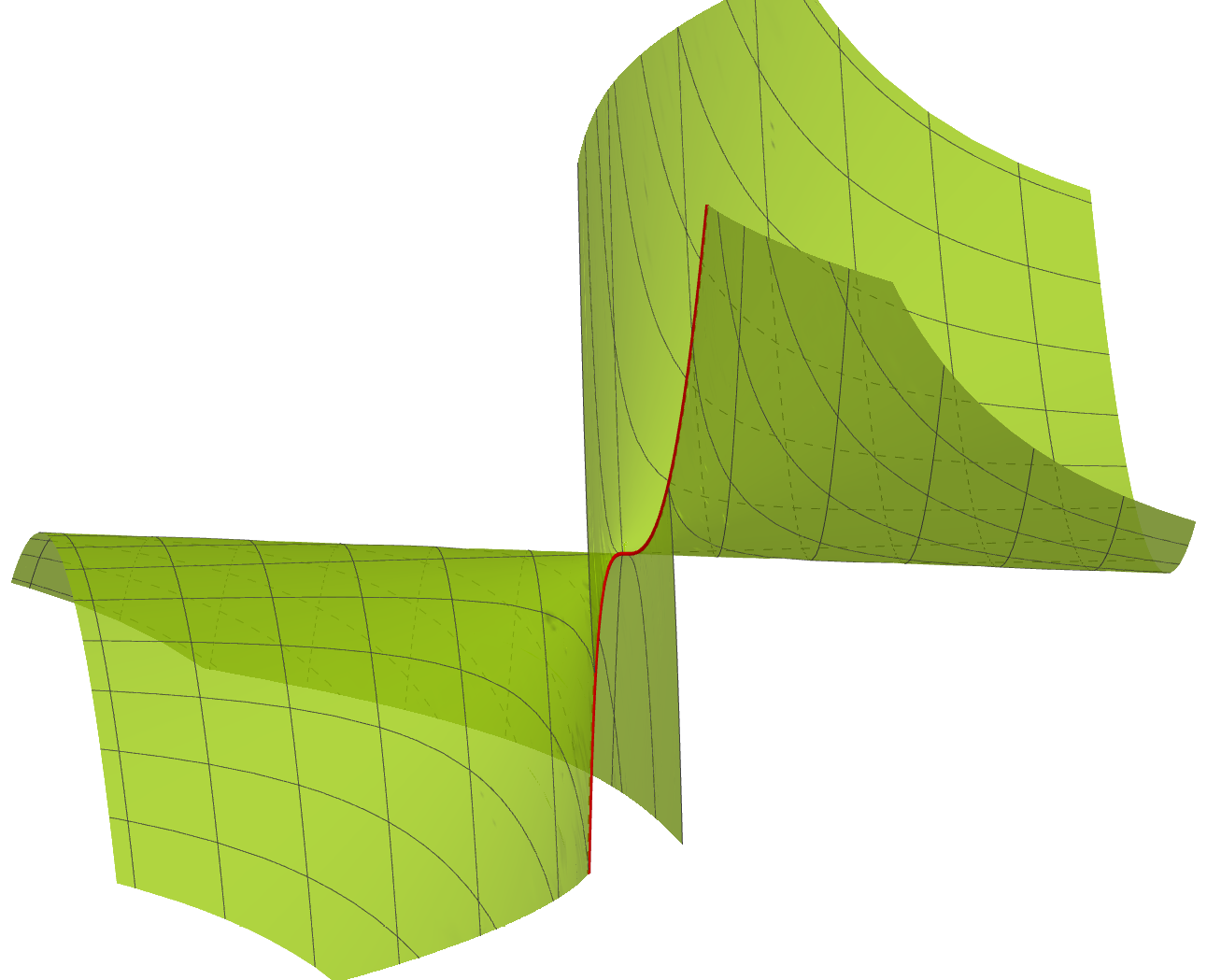}
\caption{Two partial views of the intersection of the $1$ and $2$-dimensional orbits in Einstein universe with Mink$(Y^4)$. \textbf{Red:} Part of the $1$-dimensional orbit in Minkowski patch. \textbf{Green:} Part of the $2$-dimensional orbit in Minkowski patch.}
\label{fig_1}
\end{figure}

The $1$-dimensional orbit is lightlike, homeomorphic to $\R\Pn^1$, but not a photon. The union of the $1$-dimensional orbit and the $2$-dimensional orbit is an algebraic surface, whose singular locus is precisely the $1$-dimensional orbit. It is the union of all projective lines tangent to the $1$-dimensional orbit. Figure \ref{fig_1} describes a part of the $1$ and $2$-dimensional orbits in the Minkowski patch $Mink(Y^4)$.

We will also describe the actions on Anti de-Sitter space and the generalized hyperbolic space $\Hn^{2,2}$:

\begin{theorem}\label{th.dehors}
  The orbits of PSL$(2,\R)$ in the Anti de-sitter component $\AdS^{1,3}$ are Lorentzian, i.e. of signature $(1,2)$. They are the leaves of a codimension $1$ foliation. In addition, PSL$(2,\R)$ induces three types of orbits in $\Hn^{2,2}$: a $2$-dimensional spacelike orbit (of signature $(2,0)$) homeomorphic to the hyperbolic plane $\Hn^2$, a $2$-dimensional Lorentzian orbit (i.e., of signature $(1,1)$) homeomorphic to the de-Sitter space $\dS^{1,1}$, and four kinds of $3$-dimensional orbits where the action is free:
  \begin{itemize}
  \item one-parameter family of orbits of signature $(2,1)$ consisting of elements with four distinct non-real roots,
  \item one-parameter family of Lorentzian (i.e. of signature $(1,2)$) orbits consisting of elements with four distinct real roots,
  \item two orbits of signature $(1,1,1)$,
  \item one-parameter family of Lorentzian (i.e. of signature $(1,2)$) orbits consisting of elements with two distinct real roots, and a complex root $z$ in $\mathbb{H}^2$ making an angle $\theta$ smaller than $5\pi/6$ with the two real roots.
  \end{itemize}
\end{theorem}

\begin{remark}
F. Fillastre indicated to us an alternative description for the last case stated in Theorem \ref{th.dehors}: these orbits correspond to polynomials whose roots in $\mathbb{CP}^1$ are ideal vertexes of regular ideal tetraedra in $\Hn^3$.
\end{remark}

\section{Proofs of the Theorems}

Let $f$ be an element in $\V$. We consider it as a polynomial function from $\C^2$ into $\C$. Actually, by specifying $Y=1$, we consider
$f$ as a polynomial of degree at most $4$. Such a polynomial is determined, up to a scalar, by its roots $z_1$, $z_2$, $z_3$, $z_4$ in $\C\Pn^1$
(some of these roots can be $\infty$ if $f$ can be divided by $Y$).
It provides a natural identification between $\Pn(\V)$ and the set $\widehat{\C\Pn}_4^1$ made of $4$-tuples (up to permutation) $(z_1, z_2, z_3, z_4)$ of $\C\Pn^1$ such that
if some $z_i$ is not in $\mathbb{R}\Pn^1$, then its conjugate $\bar{z}_i$ is one of the $z_j$'s.
This identification is PSL$(2,\R)$-equivariant, where the action of PSL$(2,\R)$ on $\widehat{\C\Pn}_4^1$ is simply the one induced
by the diagonal action on $(\C\Pn^1)^4$.

Actually, the complement of $\R\Pn^1$ in $\C\Pn^1$ is the union of the upper half-plane model $\mathbb{H}^2$ of the hyperbolic plane, and
the lower half-plane. We can represent every element of  $\widehat{\C\Pn}_4^1$ by a $4$-tuple (up to permutation) $(z_1, z_2, z_3, z_4)$ such that:

-- either every $z_i$ lies in $\R\Pn^1$,

-- or $z_1$, $z_2$ lies in $\R\Pn^1$, $z_3$ lies in $\mathbb{H}^2$ and $z_4 = \bar{z}_3$,

-- or $z_1$, $z_2$ lies in $\mathbb{H}^2$ and $z_3 = \bar{z}_1$, $z_4 = \bar{z}_2$.

Theorems \ref{thm1} and \ref{th.dehors} will follow from the following proposition:

\begin{proposition}\label{prop}
  Let $[f]$ be an element of $\Pn(\V)$. Then:
  \begin{itemize}
    \item it lies in $\Ein^{1,2}$ if and only if it has a root of multiplicity at least $3$, or two distinct real roots $z_1$, $z_2$, and two complex roots $z_3$, $z_4 = \bar{z}_3$, with $z_3$ in $\mathbb{H}^2$ and such that the angle at $z_3$ between the hyperbolic geodesic rays $[z_3, z_1)$ and $[z_3, z_2)$ is $5\pi/6$,
    \item it lies in $\AdS^{1,3}$ if and only it has two distinct real roots $z_1$, $z_2$, and two complex roots $z_3$, $z_4 = \bar{z}_3$, with $z_3$ in $\mathbb{H}^2$ and such that the angle at $z_3$ between the hyperbolic geodesic rays $[z_3, z_1)$ and $[z_3, z_2)$ is $> 5\pi/6$,
    \item it lies in $\mathbb{H}^{2,2}$ if and only if it has no real roots, or four distinct real roots, or a root of multiplicity exactly $2$, or it has two distinct real roots $z_1$, $z_2$, and two complex roots $z_3$, $z_4 = \bar{z}_3$, with $z_3$ in $\mathbb{H}^2$ and such that the angle at $z_3$ between the hyperbolic geodesic rays $[z_3, z_1)$ and $[z_3, z_2)$ is $< 5\pi/6$.
  \end{itemize}
\end{proposition}

\textbf{Proof of Proposition \ref{prop}.}
  Assume that $f$ has no real root. Hence we are in the situation where $z_1$, $z_2$ lie in $\mathbb{H}^2$ and $z_3 = \bar{z}_1$, $z_4 = \bar{z}_2$.
  By applying a suitable element of PSL$(2, \mathbb{R})$, we can assume $z_1=i$, and $z_2= r i$ for some $r>0$. In other words,
  $f$ is in the PSL$(2, \mathbb{R})$-orbit of $(X^2+Y^2)(X^2 + r^2Y^2)$. The value of $\mathfrak{q}$ on this polynomial is $2\times1\times r^2 + \frac16(1+r^2)^2 >0$,   hence $[f]$ lies in $\mathbb{H}^{2,2}$.

  Hence we can assume that $f$ admits at least one root in $\R\Pn^1$, and by applying a suitable element of PSL$(2, \mathbb{R})$, one can assume
  that this root is $\infty$, i.e. that $f$ is a multiple of $Y$.

  We first consider the case where this real root has multiplicity at least $2$:
  $$f = Y^2(aX^2 + bXY + cY^2)$$
  Then, $\mathfrak{q}(f) = \frac16a^2$: it follows that if $f$ has a root of multiplicity at least
  $3$, it lies in $\Ein^{1,2}$, and if it has a real root of mulitplicity $2$, it lies in $\mathbb{H}^{2,2}$.

  We assume from now that the real roots of $f$ have multiplicity $1$. Assume that all roots are real. Up to PSL$(2, \mathbb{R})$, one can assume
  that these roots are $0$, $1$, $r$ and $\infty$ with $0 < r < 1$.
  $$f(X,Y) = XY(X-Y)(X-rY) = X^3Y - (r+1)X^2Y^2 + rXY^3$$
  Then, $\mathfrak{q}(f) = -\frac12r + \frac16(r+1)^2 = \frac16(r^2-r+1) > 0$. Therefore $f$ lies in $\mathbb{H}^{2,2}$ once more.

The only remaining case is the case where $f$ has two distinct real roots, and two complex conjugate roots $z$, $\bar{z}$ with $z \in \mathbb{H}^{2}$.
Up to PSL$(2,\mathbb{R})$, one can assume that the real roots are $0$, $\infty$, hence:
$$f(X,Y) = XY(X-zY)(X-\bar{z}Y)=XY(X^2 - 2|z|\cos\theta XY + |z|^2Y^2)$$
where $z = |z|e^{i\theta}$. Then:
$$  \mathfrak{q}(f) = \frac{2|z|^2}3(\cos^2\theta - \frac34)$$
Hence $f$ lies in $\Ein^{1,2}$ if and only if $\theta = \pi/6$ or $5\pi/6$. %(observe that in this case $f(X,Y) = XY(X^2+ \sqrt{3}XY + Y^2)$).
 The proposition follows easily. \hfill$\square$

\begin{remark}
In order to determine the signature of the orbits induced by PSL$(2,\R)$ in $\Pn(\V)$, we consider the tangent vectors induced by the action of $1$-parameter subgroups of PSL$(2,\R)$. We denote by $E$, $P$ and $H$, the $1$-parameter elliptic, parabolic and hyperbolic subgroups stabilizing $i$, $\infty$ and $\{0,\infty\}$, respectively.
\end{remark}

\textbf{Proof of Theorem \ref{thm1}.} It follows from Proposition \ref{prop} that there are precisely three PSL$(2,\mathbb{R})$-orbits in
$\Ein^{1,2}$:

-- one orbit $\mathcal N$ comprising polynomials with a root of multiplicity $4$, i.e. of the form $[(sY - tX)^4]$ with $s,t \in \mathbb{R}$. It is clearly
$1$-dimensional, and equivariantly homeomorphic to $\R\Pn^1$ with the usual projective action of PSL$(2,\mathbb{R})$.
Since $\frac{d}{dt}|_{t=0}(Y-tX)^4 = -4XY^3$ is a $\mathfrak{q}$-null vector, this orbit is lightlike (but cannot be a photon since the action is irreducible),

-- one orbit $\mathcal L$ comprising polynomials with a real root of multiplicity $3$, and another real root. These are the polynomials of
the form $[(sY - tX)^3(s'Y - t'X)]$ with $s, t, s', t' \in \mathbb{R}$. It is $2$-dimensional, and it is easy to see that it is the union of the projective lines tangent to $\mathcal N$. The vectors tangent to $\mathcal{L}$ induced by the $1$-parameter subgroups $P$ and $E$ at $[XY^3]\in \mathcal{L}$ are $v_P=-Y^4$ and $v_E=3X^2Y^2+Y^4$. Obviously, $v_P$ is orthogonal to $v_E$ and $v_E+v_P$ is spacelike. Hence $\mathcal{L}$ is of signature $(0,1,1)$.

-- one open orbit comprising polynomials admitting two distinct real roots and a root $z$ in $\mathbb{H}^2$ making an angle $5\pi/6$ with the two real
roots in $\partial\mathbb{H}^2$. The stabilizers of points in this orbit are trivial since an isometry of $\mathbb{H}^2$ preserving a point in $\mathbb{H}^2$ and one point in $\partial\mathbb{H}^2$ is necessarily the identity. \hfill$\square$

\textbf{Proof of Theorem \ref{th.dehors}.}
According to Proposition \ref{prop}, the polynomials in $\AdS^{1,3}$ have two distinct real roots, and a complex root $z$ in $\mathbb{H}^2$ making an angle $\theta$ greater than $5\pi/6$ with the two real
roots. It follows that the action in $\AdS^{1,3}$ is free, and that the orbits are the level sets of the function $\theta$. Suppose that $M$ is a PSL$(2,\R)$-orbit in $\AdS^{1,3}$. There exists $r\in \R$ such that $[f]=[Y(X^2+Y^2)(X-rY)]\in M$. The orbit induced by the $1$-parameter elliptic subgroup $E$ at $[f]$ is 
$$\gamma(t)=\big[(X^2+Y^2)\big((\sin t\cos t-r \sin^2t)X^2-(\sin t\cos t+r \cos^2t)Y^2+(\cos^2t-\sin^2t+2r \sin t\cos t)XY\big)\big] .$$
Then $\mathfrak{q}(\frac{d\gamma}{dt}|_{t=0})=-2-2r^2<0$. This implies, as for any submanifold of a Lorentzian manifold admitting a timelike vector, that $M$ is Lorentzian, i.e., of signature $(1,2)$.

The case of $\mathbb{H}^{2,2}$ is the richest one. According to Proposition \ref{prop} there are four cases to consider:

\begin{itemize}
  \item \emph{No real roots.} Then $f$ has two complex roots $z_1$, $z_2$ in $\mathbb{H}^2$ (and their conjugates). It corresponds to two orbits:
  one orbit corresponding to the case $z_1 = z_2$: it is spacelike and has dimension $2$. It is the only maximal PSL$(2,\mathbb{R})$-invariant surface in $\mathbb{H}^{2,2}$ described in \cite[Section $5.3$]{Col}. The case $z_1 \neq z_2$ provides a one-parameter family of $3$-dimensional orbits on which the action is free (the parameter being the hyperbolic distance between $z_1$ and $z_2$). One may assume that $z_1=i$ and $z_2=r i$ for some $r>0$. Denote by $M$ the orbit induced by PSL$(2,\R)$ at $[f]=[(X^2+Y^2)(X^2+r^2Y^2)]$. The vectors tangent to $M$ at $[f]$ induced by the $1$-parameter subgroups $H$, $P$ and $E$ are:

\begin{align*}
v_H=-4X^4+&4r^2Y^4,\;\;\;\;\;\;\;\;\;\;\;v_P=-4X^3Y-2(r^2+1)XY^3,\\
&v_E=2(r^2-1)X^3Y+2(r^2-1)XY^3,
\end{align*}

respectively. The timelike vector $v_H$ is orthogonal to both $v_P$ and $v_E$. It is easy to see that the $2$-plane generated by $\{v_P,v_E\}$ is of signature $(1,1)$. Therefore, the tangent space $T_{[f]}M$ is of signature $(2,1)$.
  \item \emph{Four distinct real roots.} 
This case provides a one-parameter family of $3$-dimensional orbits on which the action is free - the parameter being the cross-ratio between the roots in $\R\Pn^1$. Denote by $M$ the PSL$(2,\R)$-orbit at $[f]=[XY(X-Y)(X-rY)]$ (here as explained in the proof of Proposition \ref{prop}, we can restrict ourselves to the case $0 < r < 1$). The vectors tangent to $M$ at $[f]$ induced by the $1$-parameter subgroups $H$, $P$, and $E$ are:
\begin{align*}
&v_H=-rY^4+2(r+1)XY^3-3X^2Y^2,\;\;\;\;\;\;\;v_P=-2X^3Y+2rXY^3,\\
&\;\;v_E=X^4-rY^4+3(r-1)X^2Y^2+2(r+1)XY^3-2(r+1)X^3Y,
\end{align*}
respectively. A vector $x=av_H+bv_P+cv_E$ is orthogonal to $v_P$ if and only if $2ra+b(r+1)+c(r+1)^2=0$. Set $a=\big(b(r+1)+c(r+1)^2\big)/-2r$ in
$$
\mathfrak{q}(x)=2ra^2+\frac{3}{2}b^2+\big(\frac{7}{2}(r^2+1)-r\big)c^2+2(r+1)ab+2(r+1)^2+ac(2r^2-r+5).
$$ 
Consider $\mathfrak{q}(x)=0$ as a quadratic polynomial $F$ in $b$. Since $0<r<1$, the discriminant of $F$ is non-negative and it is positive when $c\neq 0$. Thus, the intersection of the orthogonal complement of the spacelike vector $v_P$ with the tangent space $T_{[f]}M$ is a $2$-plane of signature $(1,1)$. This implies that $M$ is Lorentzian, i.e., of signature $(1,2)$.
  \item \emph{A root of multiplicity $2$.} Observe that if there is a non-real root of multiplicity $2$, when we are in the first "no real root" case. Hence we consider here only the case where the root of multiplicity $2$ lies in $\R\Pn^1$. Then, we have three subcases to consider:
\begin{itemize}
\item[--] two distinct real roots of multiplicity $2$: The orbit induced at $X^2Y^2$ is the image of the PSL$(2,\R)$-equivariant map
 $$\dS^{1,1}\subset \Pn(\R_2[X,Y])\longrightarrow \Hn^{2,2},\;\;\;\;\;\;[L]\mapsto [L^2],$$
 where $\R_2[X,Y]$ is the vector space of homogeneous polynomials of degree $2$ in two variables $X$ and $Y$, endowed with discriminant as a PSL$(2,\R)$-invariant bilinear form of signature $(1,2)$ \cite[Section $5.3$]{Col}. (Here, $L$ is the projective class of a Lorentzian bilinear form on $\R^2$). The vectors tangent to the orbit at $X^2Y^2$ induced by the $1$-parameter subgroups $P$ and $E$ are $v_P=-2XY^3$ and $v_E=2X^3Y-2XY^3$, respectively. It is easy to see that the $2$-plane generated by $\{v_p,v_E\}$ is of signature $(1,1)$. Hence, the orbit induced at $X^2Y^2$ is Lorentzian.
\item[--] three distinct real roots, one of them being of multiplicity $2$: Denote by $M$ the orbit induced by PSL$(2,\R)$ at $[f]=[XY^2(X-Y)]$. The vectors tangent to $M$ at $[f]$ induced by the $1$-parameter subgroups $H$, $P$ and $E$ are: 
$$
v_H=-2XY^3,\;\;\;\;\;v_P=Y^4-2XY^3,\;\;\;\;\;v_E=Y^4-X^4-2X^2Y^2+X^3Y-XY^3,
$$
respectively. Obviously, the lightlike vector $v_H+v_P$ is orthogonal to $T_{[f]}M$. Therefore, the restriction of the metric on $T_{[f]}M$ is degenerate. It is easy to see that the quotient of $T_{[f]}M$  by the action of the isotropic line $\R(v_H+v_P)$ is of signature $(1,1)$. Thus, $M$ is of signature $(1,1,1)$.
\item[--] one real root of multiplicity $2$, and one root in $\mathbb{H}^2$: Denote by $M$ the orbit induced by PSL$(2,\R)$ at $[f]=[Y^2(X^2+Y^2)]$. The vectors tangent to $M$ at $[f]$ induced by the $1$-parameter subgroups $H$, $P$ and $E$ are
 $$
v_H=4Y^4,\;\;\;\;\;\;v_P=-2XY^3,\;\;\;\;\;\;v_E=2X^3Y+2XY^3,$$
 respectively. Obviously, the lightlike vector $v_H$ is orthogonal $T_{[f]}M$. Therefore, the restriction of the metric on $T_{[f]}M$ is degenerate. It is easy to see that the quotient of $T_{[f]}M$  by the action of the isotropic line $\R(v_H)$ is of signature $(1,1)$. Thus $M$ is of signature $(1,1,1)$. 
\end{itemize}
  \item \emph{Two distinct real roots, and a complex root $z$ in $\mathbb{H}^2$ making an angle $\theta$ smaller than $5\pi/6$ with the two real roots.} Denote by $M$ the orbit induced by PSL$(2,\R)$ at $[f]=[Y(X^2+Y^2)(X-rY)]$. The vectors tangent to $M$ at $[f]$ induced by the $1$-parameter subgroups $H$, $P$ and $E$ are: 
\begin{align*}
v_H=-4rY^4-2X^3&Y+2XY^3,\;\;\;\;\;\;\;\;\;v_P=-3X^2Y^2+2r XY^3-Y^4,\\
&v_E=X^4-Y^4-2rX^3Y-2rXY^3,
\end{align*}

respectively. The following set of vectors is an orthogonal basis for $T_{[f]}M$ where the first vector is timelike and the two others are spacelike.
$$\{(7r+3r^3)v_H+(6-2r^2)v_P+(5+r^2)v_E,4v_P+ v_E,v_H\}.$$
Therefore, $M$ is Lorentzian, i.e., of signature $(1,2)$. \hfill $\square$
  \end{itemize}

\end{document}